\documentclass[10pt,a4paper,twoside,leqno]{article}
\usepackage{amsfonts,amssymb}

\usepackage{amscd}

\newcommand{\Stir}[2]{\left\{#1 \atop #2\right\}}%{\left\{\begin{array}{c}#1\\#2\end{array}\right\}}
\newcommand{\stir}[2]{\left[#1 \atop #2\right]}

\vspace{0.5cm}

 4 \font\ebf=cmbx8
\font\erm=cmr8 

\usepackage{graphicx}

\begin{document}

\thispagestyle{empty}

\noindent {\bf On inversion formulas  and
Fibonomial  coefficients}

\vspace{0.2cm} \noindent {\erm  A. Krzysztof Kwa\'sniewski}

\noindent {\erm member of the Institute of Combinatorics and its
Applications}

\noindent {\erm Bialystok University (*), Faculty of Physics }\\
\noindent {\erm PL - 15 - 424 Bialystok, ul. Lipowa 41,  Poland}\\
\noindent {\erm kwandr@gmail.com }

\vspace{0.2cm} \noindent {\erm  Ewa Krot-Sieniawska}

\noindent {\erm Bia{\l}ystok University, Institute of Computer Science}\\
\noindent {\erm PL - 15 - 887  Bialystok, ul. Sosnowa 64,  Poland}\\
\noindent {\erm ewakrot@wp.pl}

\noindent {\erm (*) former: Warsaw University Division }\\
\vspace{0.1cm}

\noindent {\erm FECS'08: The 2008 International Conference on
Frontiers in
 Education: WORLDCOMP'08}
\vspace{0.1cm}

\noindent {\ebf Summary}

\noindent {\small A research problem for undergraduates and
graduates is being posed as a cap for the prior antecedent regular
discrete mathematics exercises. [Here cap is not necessarily
CAP=Competitive Access Provider, though nevertheless ...] The
object of the cap problem of final interest i.e. array of
fibonomial coefficients and the issue of its combinatorial meaning
is to be found in A.K.Kwa\'sniewski's source papers. The cap
problem  number seven - still opened for students has been  placed
on Mathemagics page of the first author
[http://ii.uwb.edu.pl/akk/dydaktyka/dyskr/dyskretna.htm]. The
indicatory references are to point at a part of the vast domain of
the foundations of computer science in ArXiv affiliation
 noted as CO.cs.DM. The presentation has been verified in a tutor system of communication
 with a couple of intelligent students. The result is top secret.Temporarily.
 [Contact: Wikipedia; Theory of cognitive development].

\vspace{0.2cm}

\noindent MCS numbers: 05A19 , 11B39, 15A09 \vspace{0.2cm}


\vspace{0.2cm}

\noindent presented at the Gian-Carlo Rota Polish Seminar\\
$http://ii.uwb.edu.pl/akk/sem/sem\_rota.htm$

\vspace{0.2cm}

\subsection*{1. In the  realm of knownness. Inversion formulas. }

\textbf{Ex.1} Prove that
\begin{equation}\label{1}
\sum_{k\geq 0}{n \choose k}{k \choose l}  (-1)^{k-l}=\delta_{nl}
\end{equation}
HINT: For that to do  use the following
$$(x+1)^n=\sum_{k\geq 0}{n\choose k}x^k\;\;\Longleftrightarrow\;\; x^n=\sum_{k\geq 0}{n\choose k}(x-1)^k\;\;\Longleftrightarrow$$
$$\Longleftrightarrow\;\; x^n=\sum_{k\geq 0}\sum_{l\geq 0}{n\choose k}{k\choose l}x^l(-1)^{k-l}\;\;\;\Longleftrightarrow\;\;\;(\ref{1}).$$

\vspace{1mm}

\noindent \textbf{Ex.2} Show that
\begin{equation}\label{2}
\sum_{k=0}^n\Stir{n}{k}\stir{k}{l}(-1)^{n-k}=\delta_{nl},\;\;\;\;n\geq
0
\end{equation}
HINT: Use the combinatorial interpretation of Stirling numbers' of the I and the II kind,
($\stir{n}{k}$ and $\Stir{n}{k}$, respectively). Also note:
$$x^{\overline{n}}=\sum_{k=0}^n\stir{n}{k}x^k,\;\;\;x^n=\sum_{k=0}^n\Stir{n}{k}x^{\underline{k}},\;\;\;x^{\underline{n}}=(-1)^n(-x)^{\overline{n}}.$$
Then
$$x^k=\sum_{l=0}^k\Stir{k}{l}x^{\underline{l}}\;\;\;\wedge\;\;\;x^{\underline{l}}=(-1)^i(-x)^{\overline{l}},\;l\geq 0\;\;\Longrightarrow\;\;(\ref{2}).$$

\vspace{1mm} \noindent \textbf{Ex.3} Prove that
$$x^n=\sum_{k=0}^{n}\Stir{n}{k}(-1)^{n-k}x^{\overline{k}},\;\;\;\;x^{\underline{n}}=\sum_{k=0}^n\stir{n}{k}(-1)^{n-k}x^k,\;\;\;n\geq 0.$$
HINT: Use Exercise 2.
 \vspace{1mm}

\noindent \textbf{Ex.4}  $\left({n\choose
k}_q\right)^{-1}=?,\;\;\;$ for ${n\choose
k}_q=\frac{n_q!}{k_q!(n-k)_q!}$, $n_q\equiv \frac{1-q^n}{1-q}$,
$n,k\leq 0$.\\
This problem is solved. Just contact pp.70\&106 in \cite{1} (in
polish) and note that ${n\choose k}_q$ counts objects from the
$L(n,q)$ lattice. Then one has
$$\left({n\choose k}_q\right)^{-1}=\left({n\choose k}_q (-1)^{n-k} q^{{n-k\choose 2}}\right).$$

\vspace{1mm}

\noindent \textbf{Ex.5} Find the numbers
$\left(C_{n,\,k}\right)^{-1}$ for the sequence $C_{n,\,k}$ being
the unique solution of the recurrence relation, \cite{3}:
\begin{equation}
C_{n+1,\,k}=C_{n,\,k-1}+2C_{n,\,k}+C_{n,\,k+1}
\end{equation}
$$C_{0,\,0}=1,\;\;C_{k,\,0}=0=C_{n,\,n+k},\;\;\;\;n,k> 0.$$
For $n,k>0$ one has (\cite{3}): $C_{n,\,k}={2n\choose
n-k}\frac{k}{n}$.\\
HINT: From  the above one has: $C^{-1}_{0,\,0}=1$,
$C^{-1}_{k,\,0}=0$, $C^{-1}_{n,\,n+k}=0$ for $n,k>0$. Then one
uses $\sum_{k\geq 0}{2n\choose k}{k\choose
l}(-1)^{k-l}=\delta_{2n,\,l}$, to get
$$\sum_{s\leq n}{2n\choose s}\frac{n-s}{n}{s\choose l}\frac{n}{n-s}(-1)^{s-l}=\delta_{nl},\;\;\;n,l>0.$$

\vspace{1mm}

\noindent \textbf{Ex.6}  One can define the Charlier polynomials ,
as orthogonal polynomial sequence ( see \cite{4}, the formula
1.13, $a=-1$):
$$P_n(x)=\frac{1}{n!}\sum_{k=0}^n{n\choose k}x^{\underline{k}}.$$
Let  $x^n={\displaystyle\sum_{k\geq 0}}C_{n,k}P_k(x)$.  Find thew
numbers:
$C_{n,k}$ and $(C_{n,k})^{-1}$.\\
HINT: Use Exercise 2 and
$x^{\underline{n}}={\displaystyle\sum_{k=0}^n
}\stir{n}{k}(-1)^{n-k}x^k$.

\vspace{1mm}

\subsection*{2. Specificity beyond  the  realm  of  knownness?}}

\vspace{1mm}

\noindent \textbf{Ex.7} Discover the inversion formula i.e. the array elements   $\left({n\choose
k}_F\right)^{-1}\;\;$ for ${n\choose k}_F$ being the so called
fibonomial coefficients, i.e.
$${n\choose k}_F=\frac{n_F!}{k_F!(n-k)_F!},$$
for $n_F=F_n$ being the $n$-th Fibonacci number, ($n,k>0$). [9,5,10,11,12]\\
POSSIBLE SOLUTION: Let us consider the incidence algebra $I(\Pi)$
of the Fibonacci cobweb poset $\Pi$ and the standard reduced
incidence algebra $R(\Pi)$. These were recognized-discovered [Plato's attitude ? ] in \cite{5,6,7} and investigated there
[L.E.J. Brouwer constructivism attitude and constructivism: \\
$http://en.wikipedia.org/wiki/Constructivism\_(learning\_theory)$]\\  .
Let  $f:\Pi\times\Pi\rightarrow\mathbb{R}$ be defined as follows
$$f(x,y)=f(k,n)=\left\{\begin{array}{lll}
  {n\choose k}_F &  & k\leq n \\
  & &  \\
  0 &  & k>n \\
\end{array}\right.$$
for $x,y\in \Pi$, such that the segment $[x,y]=\{z\in\Pi:\;x\leq
z\leq y\}$ is of type $(k,n)$, i.e. $r(x)=k$, $r(y)=n$. It is
obvious that $f\in R(\Pi)$ (and of course $f\in I(\Pi)$). Then for
$g=f^{-1}$ being inverse of $f$ in $R(\Pi)$ (also in $I(\Pi)$) one
has following \cite{7}
$$
(f\ast g)(k,n)=\sum_{k\leq l\leq n}F_lg(k,l)f(l,n)=\delta_{nk},$$
i.e.
\begin{equation}
\sum_{k\leq l\leq n}F_l{n\choose l}_Fg(k,l)=\delta_{nk}.
\end{equation}
Hence in order  to discover the magic formula for  $\left({n\choose k}_F\right)^{-1}\;\;$ one has to find out
an explicit formula for $g=f^{-1}$. Right?  Maybe it can be recovered
using of the standard formula for an inverse element in $I(\Pi)$ ?,
(see for example \cite{1,8}).

%\noindent \textbf{Acknowledgements}

\end{document}